\documentclass[a4paper,11pt]{amsart}
\usepackage{hyperref,latexsym}
\usepackage{enumerate}

\theoremstyle{plain}
\newtheorem{theorem}{Theorem}[section]

\newtheorem{corollary}[theorem]{Corollary}

\theoremstyle{definition}
\newtheorem{definition}[theorem]{Definition}
\newtheorem{example}{Example}
\theoremstyle{remark}
\newtheorem{remark}{Remark}

\newcommand{\arctanh}{\rm arctanh}

\begin{document}

\title[Effective bounds in E.Hopf rigidity]
      {Effective bounds in E.Hopf rigidity for billiards and geodesic flows}

\date{20 April 2012}
\author{Misha Bialy}
\address{School of Mathematical Sciences, Raymond and Beverly Sackler Faculty of Exact Sciences, Tel Aviv University,
Israel} \email{bialy@post.tau.ac.il}
\thanks{Partially supported by ISF grant 128/10}


\begin{abstract}
In this paper we show that in some cases the E.Hopf rigidity
phenomenon admits quantitative interpretation. More precisely we
estimate from above the measure of the set $\mathcal{M}$ swept by
minimal orbits. These estimates are sharp, i.e. if $\mathcal{M}$
occupies the whole phase space we recover the E.Hopf rigidity. We
give these estimates in two cases: the first is the case of convex
billiards in the plane, sphere or hyperbolic plane. The second is
the case of conformally flat Riemannian metrics on a torus. It seems
to be a challenging question to understand such a quantitative
bounds for Burago-Ivanov theorem.
\end{abstract}

\maketitle

\section{Introduction and the result}

In this paper we estimate from above the measure of the set
$\mathcal{M}$ in the phase space which is occupied by minimal orbits
of a Hamiltonian system. These bounds are of obvious importance for
dynamics because all "rotational" invariant torii, as well as
Aubry-Mather sets are filled by minimal orbits.

The estimates provide the quantitative refinement of the E.Hopf
rigidity. We prove these bounds for two Hamiltonian systems. The
first system is a symplectic map of the cylinder corresponding to
the billiard ball motion inside a convex curve $\gamma$ lying on a
surface $\Sigma$ of constant curvature $0,\pm1$. The second system
is geodesic flow on a torus with conformally flat Riemannian metric.

Nowadays there are many cases and approaches where E.Hopf rigidity
phenomenon is established. It is an important problem to understand which of them
can be made quantitative. In particular, it seems to be a
challenging question if it is possible to give a quantitative
version for the Burago-Ivanov proof \cite{BI} of the E.Hopf conjecture.

Throughout the paper we denote by $\Omega$ the phase space of the
Hamiltonian system in question. For the billiard in a convex domain
bounded by closed curve $\gamma$, the phase space $\Omega$ is a
cylinder: $\Omega=\gamma\times(-1,1)$ equipped with the standard
symplectic form $dx\wedge d(\cos\varphi)$ giving the invariant
measure $d\mu=\sin\varphi dx d\varphi$. Here and later the billiard
map will be denoted by $T$, $x$ denotes arclength parameter on
$\gamma$ and $\varphi$ is an inward angle. As for geodesic flow on
the torus the phase space $\Omega$ is a unit tangent bundle
$\Omega=T_1\mathbf{T}^n$ equipped with the Liouville measure.

I shall use the following definition in this paper:

\begin{definition}A geodesic will be called $m-$geodesic if it has no
conjugate points. A billiard configuration $\{x_n\}$ will be called
$m$-configuration if the second variation is negative definite
between any two of end points. The corresponding orbits in the phase
space will be called $m-$orbits.
\end{definition}

Couple of remarks explain the definition. By Morse theory, for a
geodesic to be without conjugate points is equivalent to have second
variation positive definite between any two points. For billiards
any discrete Jacobi field along every $m-$configuration vanishes not
more than once and moreover change sign not more than once (see
\cite{B1} and \cite{E}).

We shall denote by $\mathcal{M}\subseteq\Omega$ the invariant subset
of the phase space consisting of all $m-$orbits. Then it follows
that $\mathcal{M}$ is a closed set (see \cite{MMS} for the discrete
case).

We shall introduce the notation for the portion of the phase space
occupied by the set
$$\Delta=\Omega\setminus\mathcal{M},\ \delta=\mu(\Delta)/\mu(\Omega),$$
where $\mu(\Omega)$ is the total measure of the phase space. Notice
that the total measure equals $2P$ for the case of billiards (here
and later $P$ denotes the length of the boundary curve $\gamma$ and
$A$ the area bounded by $\gamma$) and equals $\omega_{n-1}
Vol_g(\mathbf{T}^n)$ for a Riemannian metric $g$ on the torus (here
and below $\omega_{n-1}$ is the volume of the standard unit sphere
$\mathbf{S}^{n-1} \subset\mathbf{R}^n$. So by the definition,
$\delta\in[0,1]$ is dimensionless constant and the case $\delta=0$
is the case when all the orbits are $m-$orbits, which corresponds to
the rigidity case. The purpose of this paper is to estimate $\delta$
from below.

 We formulate first the bounds for the case of billiards:
\begin{theorem}
The following estimates hold true for a billiard table bounded by a
simple closed strictly convex curve $\gamma$ on $\Sigma$:

1. For the Euclidean plane, $\Sigma=\mathbf{R}^2$:
\begin{equation}\label{b2}
\delta\geq \frac{\pi(P^2-4\pi A)}{4P(P+\sqrt{4\pi A})}\geq
\frac{\pi(P^2-4 \pi A) }{8 P^2},
\end{equation}
and also
\begin{equation}\label{b1}
\delta\geq \frac{(P^2-4 \pi A)k_{\min} }{8 P}.
\end{equation}
 2. For the Hemisphere of,
 $\Sigma=\mathbf{S}^2$, for a curve
$\gamma$ lying entirely in the hemisphere:
\begin{equation}\label{b3}
\delta \geq
\frac{\pi}{2\arctan(\frac{1}{k_{\min}})} \frac{P^2+A^2-4\pi
A}{P(2\pi+\sqrt{P^2+(2\pi-A)^2})},
\end{equation}

3. For the Hyperbolic plane, $\Sigma=\mathbf{H}^2$, provided the
boundary curve $\gamma$ is convex with respect to horocycles, that
is $k_{\min}> 1$:
\begin{equation}\label{b4}\delta \geq
\frac{\pi}{2\arctanh(\frac{1}{k_{\min}})} \frac{P^2-A^2-4\pi
A}{P(2\pi+\sqrt{(2\pi+A)^2-P^2})}.
\end{equation}
\end{theorem}

The following remarks are in order.
\begin{remark}
 Notice that the nominators of the bounds of the theorem contain
the isoperimetric defect and therefore $\delta=0$ implies the curve
$\gamma$ is a circle on $\Sigma$. Moreover, it follows from Bonnesen
 type inequalities (see \cite{BZ}) that for small $\delta$ the curve is
close to a circle in the sense Hausdorff distance.

I would also like to mention a somewhat related result of \cite{GK} where a
quantitative version of a theorem by Mather is given estimating the
area free from caustics inside the domain bounded by $\gamma$.
\end{remark}
\begin{remark}
The estimate (\ref{b2}) uses the method of \cite{B1} where the Hopf
rigidity for billiards was found. The bounds
(\ref{b1}),(\ref{b3}),(\ref{b4}) on $\delta$ are obtained using the
so called Mirror equation. The proof of E.Hopf rigidity for plane
billiards using Mirror equation was obtained in \cite{W} and later
in \cite{B2} for the Sphere and Hyperbolic plane. Strangely the
estimates (\ref{b2}) and (\ref{b1}) are incomparable, for some
curves (\ref{b2}) is better and for others (\ref{b1}) is better. Let
me mention that it remains unclear how to push the approach of
\cite{B1} to work for Sphere and Hyperbolic plane.
\end{remark}
\begin{remark}
Let me point out that in (\ref{b4}) for the $\mathbf{H}^2$ we need
extra assumption on $\gamma$ to have $k_{\min}>1$. For the case of
rigidity when all the orbits are $m-$orbits this assumption is redundant as it is proved in \cite{B2}.
However in general case it is not clear how to get rid of it.
\end{remark}

Proof of Theorem 1.1 is given in Section 2,3.

Let me state now the result for geodesic flow. We consider
Riemannian metric on the torus $\mathbf{T}^n=\mathbf{R}^n/\Gamma$ of
the form $g=fg_0$ where $g_0$ is standard Euclidean metric on
$\mathbf{R}^n$ and $f>0$ is a conformal factor. Hopf rigidity in
this case was proved in \cite{knauf} (and later in \cite{CF} by
another method) generalizing the original proof of E.Hopf \cite{H}
and L.Green \cite{green}. Our purpose is to make their approach
quantitative and to estimate the Liouville measure $\delta$ from
below. To do this one needs a refinement of the original Hopf
method, because a straightforward application of the method does not
lead to any estimate on $\delta$ (it is especially clear for the
case $n=2$). For the proof I used below some of the earlier ideas of
\cite{B3} on rigidity of Newton equations.

We shall split the result into two cases, $n=2$ and $n>2$.

\begin{theorem}\label{t2}1. For $n=2$, let $\psi:\mathbf{R}_+\rightarrow\mathbf{R}_+$
be any positive smooth function. Denote by
$\Psi(f)=\psi'(f)\left(\frac{4}{f}-\frac{\psi'(f)}{\psi(f)}\right)$.
Then the following estimate holds true:
$$\delta\geq\frac{\pi
\int_{\mathbf{T}^2} \Psi(f)|grad_{g_0}f|_{g_0}^2\
dVol_{g_o}}{4\|K\|_{C^0} \|\psi(f)\|_{C^0}Vol(\mathbf{T}^2,g)},
$$
where $K$ is the curvature of the metric $g$.

2. For $n>2$, for any positive function
$\psi:\mathbf{R}_+\rightarrow\mathbf{R}_+$ introduce $$\Psi(f)=
\Psi(f)=f^{\frac{n}{2}-1}\psi'(f)\left(\frac{4}{f}-\frac{\psi'(f)}{\psi(f)}\right)
+(n-2)f^{\frac{n}{2}-3}\psi(f).$$
Then the following estimate holds

$$\delta\geq\frac{(n-1)\omega_{n-1}
\int_{\mathbf{T}^n} \Psi(f)|grad_{g_0}f|_{g_0}^2\
dVol_{g_o}}{4n\|Ric\|_{C^0} \|\psi(f)\|_{C^0}Vol(\mathbf{T}^n,g)},
$$
where $Ric$ stands for the Ricci tensor of $g$.
\end{theorem}
Obviously this statement makes sense only if $\Psi$ is positive
function. It turns out to be possible for many choices of $\psi$.
\begin{corollary}
For the particular choice of $\psi(f)=f^\alpha$ we
have:

1. For $n=2$ and for every $\alpha$ in the range $0<\alpha<4$ it follows \newline
$\Psi(f)=\alpha(4-\alpha)f^{\alpha-2}$ and thus
$$\delta\geq\frac{\pi\alpha(4-\alpha)
\int_{\mathbf{T}^2} f^{\alpha-2}(f_{x_1}^2+f_{x_2}^2)
dx_1dx_2}{4\|K\|_{C^0} \|f\|_{C^0}^\alpha\int fdx_1dx_2}.
$$

2. For $n>2$ and for every $\alpha$ in the range where $(n-2)+\alpha(4-\alpha)>0$
it follows
$$\Psi(f)=((n-2)+\alpha(4-\alpha))f^{\frac{n}{2}-3+\alpha}$$ and thus
$$\delta\geq((n-2)+\alpha(4-\alpha))\frac{(n-1)\omega_{n-1}
\int_{\mathbf{T}^n} f^{\frac{n}{2}-3+\alpha}|grad_{g_0}f|_{g_0}^2\
dVol_{g_o}}{4n\|Ric\|_{C^0} \|f^\alpha\|_{C^0}Vol(\mathbf{T}^n,g)},
$$
\end{corollary}

\begin{example}
For $n=2$ and $\alpha=2$ one has
$$\delta\geq\frac{\pi
\int_{\mathbf{T}^2}(f_{x_1}^2+f_{x_2}^2) dx_1dx_2}{\|K\|_{C^0}
\|f\|_{C^0}^2\int fdx_1dx_2},
$$
As for $n>2$ and $\alpha=2$ one has:
$$\delta\geq\frac{(n+2)(n-1)\omega_{n-1}
\int_{\mathbf{T}^n} f^{\frac{n}{2}-1}|grad_{g_0}f|_{g_0}^2\
dVol_{g_o}}{4n\|Ric\|_{C^0} \|f\|^2_{C^0}Vol(\mathbf{T}^n,g)}.
$$
\end{example}
\begin{remark} For both cases $n=2$ and $n>2$ and for $\alpha \geq 2$
one gets the strongest estimate in the Corollary for $\alpha=2$,
because then the value of $\alpha(4-\alpha)$ becomes maximal.
Analogously, for the case $n>2$ the estimate of the Corollary for
$\alpha\leq 0$ is best possible for $\alpha=0$. Thus the meaningful
range for $\alpha$ in the Corollary is $\alpha\in(0,2]$, for $n=2$
and $\alpha\in[0,2]$ for $n\geq 2$. Besides these remarks the
estimates for different values of $\alpha$ seem to be incomparable.
 Let me also point out that unlike
the case $n=2$, for $n>2$ the choice $\alpha=0$ is allowed and
corresponds to the inequality considered by A.Knauf.
\end{remark}
Proofs of Theorem \ref{t2} are given in Sections 4,5.
\section*{Acknowledgements} It is a pleasure to thank Semyon Alesker, Maxim Arnold, Victor
Bangert and Wilderich Tuschmann for valuable discussions.

\section{Estimates for planar billiards.}
It follows from \cite{B1} that along any m-configuration one can
construct a positive discrete Jacobi field and then using this field
to define a bounded measurable function
$\omega:\mathcal{M}\rightarrow\mathbf{R},$ satisfying the
inequality:
\begin{equation}\label{in1}
 \omega(y, \psi)-\omega(x,
\varphi) \geq L_{11} (x, y)+ 2L_{12}(x, y)+ L_{22}(x, y).
\end{equation}
 Here and below $T:(x,cos\varphi)\mapsto(y,\cos\psi)$; $L$
denotes the distance between $\gamma(x)$ and $\gamma(y)$; $x$ is an
arclength parameter on $\gamma$ and subindexes of $L$ stand for
partial derivatives with respect to $x,y$ respectively.

Integrate against the invariant measure $\mu$ inequality (\ref{in1})
over the set $\mathcal{M}$ of all $m-$orbits. We get:
$$0\geq \int_{\mathcal{M}}(L_{11} (x, y)+ 2L_{12}(x, y)+ L_{22}(x, y))d\mu.$$
After a computation this leads to the inequality:
\begin{equation}
\int_{\mathcal{M}}\frac{(\sin\varphi+\sin\psi)^2}{L}d\mu\leq
\int_{\mathcal{M}}(k(x)\sin\varphi+k(y)\sin\psi)\ d\mu.\label{in2}
\end{equation}
The LHS of (\ref{in2}) can be estimated from below by Cauchy-Schwartz
and Santalo formulas:
$$LHS\geq
\frac{(\int_{\mathcal{M}}(\sin\varphi+\sin\psi)d\mu)^2}{\int_{\mathcal{M}}{L}d\mu}\geq
\frac{(2\int_{\mathcal{M}}\sin\varphi \
d\mu)^2}{\int_{\Omega}{L}d\mu}=
\frac{(2\int_{\mathcal{M}}\sin\varphi\ d\mu)^2}{2\pi A}.
$$

The RHS of (\ref{in2}) can be estimated:
$$
RHS=2\int_{\mathcal{M}}k(x)\sin\varphi\ d\mu\leq
2\int_{\Omega}k(x)\sin\varphi\ d\mu=2\pi^2.
$$
Therefore (\ref{in2}) gives the following
$$
\frac{2\int_{\mathcal{M}}\sin\varphi\  d\mu}{\sqrt{2\pi A}}\leq
\sqrt{2\pi^2}.$$ Therefore
$$
2\int_{\mathcal{M}}\sin\varphi\  d\mu\leq \pi\sqrt{4\pi A} ,
$$
Estimating the left hand side of the last inequality we get
$$
\pi P-4\delta P\leq 2\int_{\Omega}\sin\varphi\
d\mu-2\int_{\Delta}\sin\varphi  d\mu=
2\int_{\mathcal{M}}\sin\varphi\ d\mu\leq \pi\sqrt{4\pi A}.
$$
Thus
$$
\frac{\sqrt{4\pi A}}{P}\geq 1-\frac{4}{\pi}\delta,
$$
so that
$$\frac{\pi}{4}\left(1-\frac{\sqrt{4\pi A}}{P}
\right)\leq\delta.
$$
Then
$$\frac{\pi}{4}\left(\frac{P^2-4\pi A}{P(2P)}\right)\leq
\frac{\pi}{4}\left(\frac{P^2-4\pi A}{P(P+\sqrt{4\pi
A})}\right)\leq\delta.
$$
This proves (\ref{b2}).

 In order to prove (\ref{b1}) we use another measurable function defined on the
subset filled by $m-$orbits $$a:\mathcal{M}\rightarrow\mathbf{R},\
0<a(x,\varphi)<L(x,\varphi),$$ which is related in fact to the
function $\omega$ discussed in the proof of (\ref{b2})(see \cite{B2}). This
function satisfies the Mirror equation for any point
$(x,\varphi)\in\mathcal{M}$:
\begin{equation}
\frac{1}{a(x,\varphi)}+\frac{1}{
L(T^{-1}(x,\varphi))-a(T^{-1}(x,\varphi))}=\frac{2k(x)}{\sin{\varphi}},\label{mp}
\end{equation}
Then it follows
$$\frac{a(x,\varphi) + (L(T^{-1}(x,\varphi))-a(T^{-1}(x,\varphi))}{2}\geq\frac{\sin\varphi}{k(x)}.
$$
Integrate this inequality against the invariant measure $d\mu$ over
the set $\mathcal{M}.$ We have:
\begin{equation}\label{in3}
\frac{1}{2}\int_{\mathcal{M}}L
d\mu\geq\int_{\mathcal{M}}\frac{\sin\varphi}{k(x)}.
\end{equation}
The LHS of (\ref{in3}) can be estimated using Santalo formula:
$$\pi A=\frac{1}{2}\int_{\Omega}L d\mu \geq
\frac{1}{2}\int_{\mathcal{M}}L d\mu.
$$
And for the RHS we have using Cauchy-Schwartz
$$
\int_{\mathcal{M}}\frac{\sin\varphi}{k(x)}=\int_{\Omega}\frac{\sin\varphi}{k(x)}-
\int_{\Delta}\frac{\sin\varphi}{k(x)}\geq
$$
$$\quad\quad\geq\frac{\pi}{2}\int_0^P\frac{1}{k(x)}dx-\frac{2\delta P}{k_{\min}}\geq
\frac{\pi}{2}\cdot \frac{P^2}{2\pi}-\frac{2\delta P}{k_{\min}}.
$$
Therefore (\ref{in3}) yields:
$$\pi A\geq \frac{P^2}{4}-\frac{2\delta P}{k_{\min}},$$
which is equivalent to (\ref{b1}). This completes the proof of (\ref{b1}).

\section{Billiard on the Sphere and the Hyperbolic plane}
The Mirror equation for billiards on Hemisphere and Hyperbolic plane
is obtained in \cite{B2}. For the Hemisphere, there exists a
measurable function $$a:\mathcal{M}\rightarrow\mathbf{R},\
0<a(x,\varphi)<L(x,\varphi)$$ such that for any point
$(x,\varphi)\in\mathcal{M}$ the following holds:
\begin{equation}
\cot \left(a(x,\varphi)\right)+\cot\left
(L(T^{-1}(x,\varphi))-a(T^{-1}(x,\varphi))\right)=
\frac{2k(x)}{\sin{\varphi}}.\label{mhemi}
\end{equation}
It implies:
$$
\cot
\frac{a(x,\varphi)+L(T^{-1}(x,\varphi))-a(T^{-1}(x,\varphi))}{2}\leq\frac
{k(x)}{\sin\varphi},
$$
Equivalently
$$\frac{a(x,\varphi)+
L(T^{-1}(x,\varphi))-a(T^{-1}(x,\varphi))}{2}\geq
\arctan\left(\frac{\sin\varphi}{k(x)}\right).
$$
Integrating over $\mathcal{M}$ with respect to the invariant measure
$d \mu=\sin\varphi\ dxd\varphi$ we get: \begin{equation}\label{in4}
\int_{\mathcal{M}} L\ d\mu\ \geq\ 2\int_{\mathcal{M}}
\arctan\left(\frac{\sin\varphi}{k(x)}\right)d\mu
\end{equation}
For the LHS of (\ref{in4}) we have:
$$
\int_{\mathcal{M}} L\ d\mu\leq\int_{\Omega} L\ d\mu=2\pi A.
$$
As for the RHS of (\ref{in4}) we compute and use Gauss-Bonnet to
get:
$$
2\int_{\mathcal{M}}
\arctan\left(\frac{\sin\varphi}{k(x)}\right)d\mu= 2\int_{\Omega}
\arctan\left(\frac{\sin\varphi}{k(x)}\right)d\mu-2\int_{\Delta}
\arctan\left(\frac{\sin\varphi}{k(x)}\right)d\mu
$$
$$\geq 2\int_0^P dx\int_0^{\pi}
\arctan\left(\frac{\sin\varphi}{k(x)}\right)\sin\varphi\
d\varphi-4P\delta\arctan\left(\frac{1}{k_{\min}}\right)=
$$
$$
=2\pi\int_0^P(\sqrt{k^2(x)+1}-k(x))dx-
4P\delta\arctan\left(\frac{1}{k_{\min}}\right)=
$$
$$
=2\pi\int_0^P\sqrt{k^2(x)+1}dx-2\pi(2\pi-A)-
4P\delta\arctan\left(\frac{1}{k_{\min}}\right).
$$
Substitute now the estimates back into (\ref{in4}):
$$
\int_0^P\sqrt{k^2(x)+1}dx\leq
2\pi+\frac{2P\delta}{\pi}\arctan\left(\frac{1}{k_{\min}}\right).
$$
But then the following two inequalities follow. The first one:
$$
\int_0^P(\sqrt{k^2(x)+1}-1)dx\leq 2\pi-P+
\frac{2P\delta}{\pi}\arctan\left(\frac{1}{k_{\min}}\right).
$$
 And the second is:
$$
\int_0^P(\sqrt{k^2(x)+1}+1)dx\leq 2\pi+P+
\frac{2P\delta}{\pi}\arctan\left(\frac{1}{k_{\min}}\right).
$$
Multiplying two of them and using Cauchy-Schwartz we get:
$$
(2\pi-A)^2=\left(\int_0^Pk(x)dx\right)^2\leq \left(2\pi+
\frac{2P\delta}{\pi}\arctan\left(\frac{1}{k_{\min}}\right)\right)^2-P^2.
$$
Therefore $$ \sqrt{P^2+(2\pi -A)^2}\leq 2\pi+
\frac{2P\delta}{\pi}\arctan\left(\frac{1}{k_{\min}}\right).
$$
Thus
$$
\frac{P^2-4\pi A+A^2}{\sqrt{P^2+(2\pi -A)^2}+2\pi}= \sqrt{P^2+(2\pi
-A)^2}-2\pi\leq
$$
$$
\leq\frac{2P\delta}{\pi}\arctan\left(\frac{1}{k_{\min}}\right).
$$
But this is exactly (\ref{b3}), so the proof for the Hemisphere is
finished.

 For the Hyperbolic plane the proof is similar. Let me sketch the
main steps. Let me remind first that for the Hyperbolic case we need an additional
requirement $k(x)\geq k_{\min}>1$. In particular this implies
$$P< k_{\min} P\leq\int_0^P k(x)dx=2\pi+A.$$

We start again with a measurable function
$$a:\mathcal{M}\rightarrow\mathbf{R},\ 0<a(x,\varphi)<L(x,\varphi)$$ such
that for any point $(x,\varphi)\in\mathcal{M}$ the Mirror equation
holds:
\begin{equation}
\coth \left(a(x,\varphi)\right)+\coth\left
(L(T^{-1}(x,\varphi))-a(T^{-1}(x,\varphi))\right)=
\frac{2k(x)}{\sin{\varphi}}.\label{mhyp}
\end{equation}
This leads to the inequality: $$ a(x,\varphi)+
L(T^{-1}(x,\varphi))-a(T^{-1}(x,\varphi))\geq 2
\arctanh\left(\frac{\sin\varphi}{k(x)}\right).
$$
Integrating over $\mathcal{M}$ we get
$$
\int_{\mathcal{M}} L\ d\mu\ \geq\ 2\int_{\mathcal{M}} \arctanh
\left(\frac{\sin\varphi}{k(x)}\right)d\mu,
$$
which leads to the inequality:
$$2\pi A\geq 2\int_{\mathcal{M}} L\ d\mu\ \geq\pi \int_0^P(k(x)-\sqrt{k^2(x)-1})dx-
4P\delta\arctanh \left(\frac{1}{k_{\min}}\right).
$$
This implies by Gauss-Bonnet:
$$
\int_0^P\sqrt{k^2(x)-1}dx\geq 2\pi-\frac{2P\delta}{\pi}\arctanh
\left(\frac{1}{k_{\min}}\right).
$$
By Cauchy Schwartz we have;
$$
\int_0^P\sqrt{k^2(x)-1}dx\leq\left(\int_0^P(
k(x)-1)dx\int_0^P(k(x)+1)dx\right)^\frac{1}{2}=
\sqrt{(A+2\pi)^2-P^2}.
$$
Thus we get:
$$
\sqrt{(A+2\pi)^2-P^2}\geq2\pi-\frac{2P\delta}{\pi}\arctanh
\left(\frac{1}{k_{\min}}\right).
$$
This completes the proof.

\section{Proof of the estimates for geodesic flows in $n=2$.}
The original E.Hopf method needs a modification in order to get
bounds on the measure $m-$geodesics. This goes as follows.

First, following E.Hopf, for every geodesic with no conjugate points
one constructs by a limiting procedure a positive solution of the
Jacobi equation and then a measurable bounded function
$\omega:\mathcal{M}\rightarrow\mathbf{R}$ which is smooth along the
orbits of the geodesic flow satisfying the Riccati equation:
\begin{equation}
\dot{\omega}+\omega^2+K=0.
\end{equation}
Here the derivative is taken in the direction of the vector of the
geodesic flow in $T_1\mathbf{T}^2$, and $K$ is the curvature of the
conformal metric $g=f(dx_1^2+dx_2^2)$. Let me remind that
$\mathcal{M}$ is a closed  subset of the phase space
$\Omega=T_1\mathbf{T}^2$ invariant under the geodesic flow.

Multiplying both sides of the equation by a positive factor
$\psi(f)$ we get:
\begin{equation}
\frac{d}{dt}(\psi\omega)-\omega\frac{d}{dt}(\psi(f))+\psi(f)\omega^2
+\psi(f) K=0,
\end{equation}
Which leads to
\begin{equation}
\frac{d}{dt}(\psi(f)\omega)- \psi^{'}(f)(f_{x_1}
\dot{x_1}+f_{x_2}\dot{x_2})\omega+\psi(f)\omega^2 +\psi(f) K=0,
\end{equation}

For $T_1\mathbf{T}^2$ we have $\dot{x_1}=\frac{1}{\sqrt
f}\cos\varphi, \dot{x_2}=\frac{1}{\sqrt f}\sin\varphi$ therefore

\begin{equation}
\frac{d}{dt}(\psi(f)\omega)- \psi^{'}(f)\left( \frac{f_{x_1}}{\sqrt
f}\cos\varphi+\frac{f_{x_2}}{\sqrt
f}\sin\varphi\right)\omega+\psi(f)\omega^2 +\psi(f) K=0,
\end{equation}
Integrate the last equation over the set $\mathcal{M}$ against the
Liouville measure $d\mu=fdx_1dx_2d\varphi$ and use its invariance
under the geodesic flow to get

\begin{equation}\label{eq}
-\int_{\mathcal{M}} \psi^{'}\left( \frac{f_{x_1}}{\sqrt
f}\cos\varphi+\frac{f_{x_2}}{\sqrt f}\sin\varphi\right)\omega
d\mu+\int_{\mathcal{M}}\psi\omega^2d\mu +\int_{\mathcal{M}}\psi
Kd\mu=0,
\end{equation}

Denote the first and the last term in equation (\ref{eq}) by $A$ and
$C$ respectively. Then, by Cauchy-Schwartz for $A$ we have:

$$A\geq
-\left(\int_{\mathcal{M}} \frac{(\psi^{'})^2}{f\psi}\left(
f_{x_1}\cos\varphi+f_{x_2}\sin\varphi\right)^2d\mu\right)^\frac{1}{2}
\left(\int_{\mathcal{M}}\psi\omega^2d\mu \right)^\frac{1}{2}\geq
$$
$$\quad -\left(\int_{\Omega} \frac{(\psi^{'})^2}{f\psi}\left(
f_{x_1}\cos\varphi+f_{x_2}\sin\varphi\right)^2d\mu\right)^\frac{1}{2}
\left(\int_{\mathcal{M}}\psi\omega^2d\mu \right)^\frac{1}{2}=
$$
$$-\left(\pi\int_{\mathbf{T}^2} \frac{(\psi^{'})^2}{\psi}(
f_{x_1}^2+f_{x_2}^2)dx_1dx_2\right)^\frac{1}{2}
\left(\int_{\mathcal{M}}\psi\omega^2d\mu \right)^\frac{1}{2}
$$
The third term $C$ can be written as follows:
$$C=\int_{\Omega}\psi Kd\mu-\int_{\Delta}\psi Kd\mu\geq
\int_{\Omega}\psi Kd\mu-\|\psi(f)\|_{C_0}\|K\|_{C_0}\mu(\Delta)=
$$
$$
=2\pi\int_{\mathbf{T}^2}\psi(f)K fdx_1dx_2-\|\psi(f)\|_{C_0}\|K\|_{C_0}\cdot\delta\cdot Vol(\mathbf{T}^2,g).
$$
Substitute the explicit expression for $K=-\frac{\Delta(\log
f)}{2f}$ and integrate by parts to get:

$$
C\geq \pi \int_{\mathbf{T}^2}\frac{\psi'(f)}{f}(
f_{x_1}^2+f_{x_2}^2)dx_1dx_2-\|\psi(f)\|_{C_0}\|K\|_{C_0}\delta\cdot Vol(\mathbf{T}^2,g).
$$
Use the estimates of the terms $A$ and $C$ in the equation (\ref{eq})
$$
-\left(\pi\int_{\mathbf{T}^2} \frac{(\psi^{'})^2}{\psi}(
f_{x_1}^2+f_{x_2}^2)dx_1dx_2\right)^\frac{1}{2}\cdot X +X^2+
$$
\begin{equation}\label{inequality}
+\pi \int_{\mathbf{T}^2}\frac{\psi'(f)}{f}(
f_{x_1}^2+f_{x_2}^2)dx_1dx_2-
\|\psi(f)\|_{C_0}\|K\|_{C_0}\cdot\delta\cdot Vol(\mathbf{T}^2,g)\leq
0,
\end{equation}
where we denoted by
$$
X=\left(\int_{\mathcal{M}}\psi\omega^2d\mu \right)^\frac{1}{2}.
$$
Next notice that (\ref{inequality}) is a quadratic inequality in $X$
and therefore the discriminant must be non-negative:

$$\pi\int_{\mathbf{T}^2} \frac{(\psi^{'})^2}{\psi}(
f_{x_1}^2+f_{x_2}^2)dx_1dx_2-4\pi
\int_{\mathbf{T}^2}\frac{\psi'(f)}{f}( f_{x_1}^2+f_{x_2}^2)dx_1dx_2+
$$
$$
\quad\quad\quad+4 \|\psi(f)\|_{C_0}\|K\|_{C_0}\cdot\delta\cdot
Vol(\mathbf{T}^2,g)\geq 0.
$$
But this is precisely the inequality which is claimed. This
completes the proof for $n=2$.

\section{Proof of the estimates for geodesic flows in $n>2$.}
In this case we modify the approach by L.Green and A.Knauf in a
similar way that we did for the case $n=2$. We start with a
measurable bounded function (see \cite{green} or \cite{chavel} for
the construction) $\omega:\mathcal{M}\rightarrow\mathbf{R}$ which
satisfies the differential inequality:
$$\frac{d}{dt}{\omega}+\frac{\omega^2}{n-1}+R\leq 0,$$
where the derivative is along the geodesic flow and $R$ is a
function $$R:\Omega\rightarrow\mathbf{R}, R(v)=Ric(v,v).$$
Multiplying both sides of the inequality by a positive factor
$\psi(f)$ we get:
\begin{equation}
\frac{d}{dt}(\psi\omega)-\psi^{'}(f)\dot{f}\omega+
\psi(f)\frac{\omega^2}{n-1} +\psi(f) R\leq 0,
\end{equation}

Integrate against the invariant measure $d\mu=f^{\frac{n}{2}}dxdo$
over the set $\mathcal{M}\subseteq \Omega$ (where $dx$, $do$ are the
standard measures on Euclidean space and on the unit sphere).
\begin{equation}\label{eq1}
-\int_{\mathcal{M}} \psi^{'}(f)\dot{f}\omega
d\mu+\int_{\mathcal{M}}\frac{\psi\omega^2}{n-1}d\mu
+\int_{\mathcal{M}}\psi Rd\mu\leq 0,
\end{equation}
We can estimate the first term $A$ and the last $C$ as follows. By
Cauchy-Schwartz inequality

$$A\geq -\left(\int_{\mathcal{M}} \frac{(\psi^{'})^2}{\psi}(\dot{f})^2d\mu\right)^\frac{1}{2}
\left(\int_{\mathcal{M}}\psi\omega^2d\mu \right)^\frac{1}{2}\geq
$$
$$-\left(\int_{\Omega} \frac{(\psi^{'})^2}{\psi}(\dot{f})^2d\mu\right)^\frac{1}{2}
\left(\int_{\mathcal{M}}\psi\omega^2d\mu \right)^\frac{1}{2}=
$$
$$-\left(\int_{\Omega}\frac{(\psi^{'})^2}{f\psi}<grad_{g_0} f,
f\dot{x}>_{g_0}^2 f^{\frac{n}{2}}dxdo\right)^\frac{1}{2}
\left(\int_{\mathcal{M}}\psi\omega^2d\mu \right)^\frac{1}{2}=
$$
$$
-\left(\frac{\omega_{n-1}}{n}\int_{\mathbf{T}^n}\frac{(\psi^{'})^2}{\psi}\|grad_{g_0}
f\|_{g_0}^2 f^{\frac{n}{2}-1}dx\right)^\frac{1}{2}
\left(\int_{\mathcal{M}}\psi\omega^2d\mu \right)^\frac{1}{2}.
$$
For the last term $C$ we have:
$$C=\int_{\Omega}\psi Rd\mu-\int_{\Delta}\psi R d\mu\geq
\int_{\Omega}\psi Rd\mu-\|\psi(f)\|_{C_0}\|R\|_{C_0}\mu(\Delta)=
$$
$$
=\frac{\omega_{n-1}}{n}\int_{\mathbf{T}^n}\psi(f)Scal(g)
f^{\frac{n}{2}}dx-\|\psi(f)\|_{C_0}\|R\|_{C_0}\cdot\delta\cdot
Vol(\mathbf{T}^n,g),
$$
where $Scal(g)$ is the Scalar curvature of $g$. Substitute the
explicit expression for $Scal$,
$$Scal(g)=(1-n)f^{-2}\Delta f+\frac{(1-n)(n-6)}{4}f^{-3}\|grad_{g_0}f\|_{g_0}^2$$
and integrate by parts the term with the Laplacian to get:

$$
C\geq
\frac{\omega_{n-1}(1-n)}{n}\left(\int_{\mathbf{T}^n}\psi(f)f^{\frac{n}{2}-2}
\Delta fdx
+\int_{\mathbf{T}^n}\frac{(n-6)}{4}\psi(f)f^{\frac{n}{2}-3}\|grad_{g_0}f\|_{g_0}^2dx\right)-
$$
$$
-\|\psi(f)\|_{C_0}\|R\|_{C_0}\cdot\delta\cdot
Vol(\mathbf{T}^n,g)\geq $$
$$\geq\frac{\omega_{n-1}(1-n)}{n}\int_{\mathbf{T}^n}
\left(-(\psi(f)f^{\frac{n}{2}-2})^{'}+
\frac{(n-6)}{4}\psi(f)f^{\frac{n}{2}-3}\right)\|grad_{g_0}f\|_{g_0}^2dx-
$$
$$
-\|\psi(f)\|_{C_0}\|R\|_{C_0}\cdot\delta\cdot
Vol(\mathbf{T}^n,g)=:\tilde{C}.$$

Substituting the estimates on $A,C$ and the notation
$X=\left(\int_{\mathcal{M}}\psi(f)\omega^2d\mu\right)^\frac{1}{2} $
into (\ref{eq1}) we get the quadratic inequality:

$$
-X\cdot\left(\frac{\omega_{n-1}}{n}\int_{\mathbf{T}^n}\frac{(\psi^{'})^2}{\psi}\|grad_{g_0}
f\|_{g_0}^2 f^{\frac{n}{2}-1}dx\right)^\frac{1}{2}
+\frac{1}{n-1}X^2+\tilde{C}\leq 0.
$$
So the discriminant $D$ must be non-negative, which leads to

$$
\frac{\omega_{n-1}}{n}\int_{\mathbf{T}^n}\frac{(\psi^{'})^2}{\psi}\|grad_{g_0}
f\|_{g_0}^2 f^{\frac{n}{2}-1}dx-\frac{4}{n-1}\tilde{C}\geq 0.
$$
Then
$$
  \frac{4}{n-1}
\|\psi(f)\|_{C_0}\|R\|_{C_0}\cdot\delta\cdot Vol(\mathbf{T}^n,g)\geq
$$
$$
-\frac{\omega_{n-1}}{n}\int_{\mathbf{T}^n}\left(\frac{(\psi^{'})^2}{\psi}
f^{\frac{n}{2}-1}-4(\psi(f)f^{\frac{n}{2}-2})^{'}+
(n-6)\psi(f)f^{\frac{n}{2}-3}\right)\|grad_{g_0}f\|_{g_0}^2dx.
$$
By the definition of $\Psi(f)$ this is equivalent to:
$$\frac{4}{n-1}
\|\psi(f)\|_{C_0}\|R\|_{C_0}\cdot\delta\cdot Vol(\mathbf{T}^n,g)\geq
\frac{\omega_{n-1}}{n}\int_{\mathbf{T}^n}\Psi(f)\|grad_{g_0}f\|_{g_0}^2dx.
$$
This proves the claim for $n>2$.

\end{document}